\documentclass{ifacconf}
\usepackage{natbib}
\usepackage{graphicx}      
\usepackage{amssymb}
\usepackage{amstext}
\usepackage{amsmath}
\usepackage{float}
\usepackage{color}

\def\sign{\mbox{sign\,}}

\newtheorem{example}{Example}

\begin{document}
\begin{frontmatter}

\title{A short survey on QPSK Costas loop mathematical models}

\author[spb,fin]{Kuznetsov N.V.}
\author[spb]{Kuznetsova O.A.}
\author[spb,ipmash]{Leonov G.A.}
\author[spb]{Yuldashev M.V.}
\author[spb]{Yuldashev R.V.}

\address[spb]{Faculty of Mathematics and Mechanics, Saint-Petersburg State University, Russia 
}
\address[fin]{Dept. of Mathematical Information Technology,University of Jyv\"{a}skyl\"{a}, Finland (e-mail: nkuznetsov239@gmail.com)}
\address[ipmash]{Institute for Problems in Mechanical Engineering of the Russian Academy of Sciences, Russia}

\begin{abstract}
The Costas loop is a modification of the phase-locked loop circuit,
which demodulates data and recovers carrier from the input signal.
The Costas loop is essentially a nonlinear control system and
its nonlinear analysis is a challenging task.
Thus, simplified mathematical models and their numerical simulation
are widely used for its analysis.
At the same time for phase-locked loop circuits there are known various
examples where the results of such simplified analysis
are differ substantially from the real behavior of the circuit.
In this survey the corresponding problems
are demonstrated and discussed for the QPSK Costas loop.
\end{abstract}

\end{frontmatter}

\section{Introduction}
The Costas loop is a classical modification of the phase-locked loop circuit (PLL),
which is a nonlinear control system designed to generate an electrical signal,
the phase of which is automatically tuned to the phase of the input signal.
The Costas loop is essentially a nonlinear control system and
its nonlinear analysis is a challenging task.
Thus, simplified mathematical models and their numerical simulation
are widely used for its analysis.
At the same time for PLL based circuits there are known various
examples where the results of such simplified analysis
are differ substantially from the real behavior of the circuit
(see corresponding discussion of gaps between mathematical control
theory, the theory of dynamical systems and the engineering practice of PLL in \citep{LeonovKYY-2015-TCAS}).
Recently such examples were discussed for the BPSK Costas loop
in \citep{BestKKLYY-2015-ACC,BestKLYY-2016}.
In this survey the corresponding problems
are revealed and discussed for the QPSK Costas loop.

\section{QPSK Costas loop operation}

\begin{figure}[h]
  \includegraphics[scale=0.4]{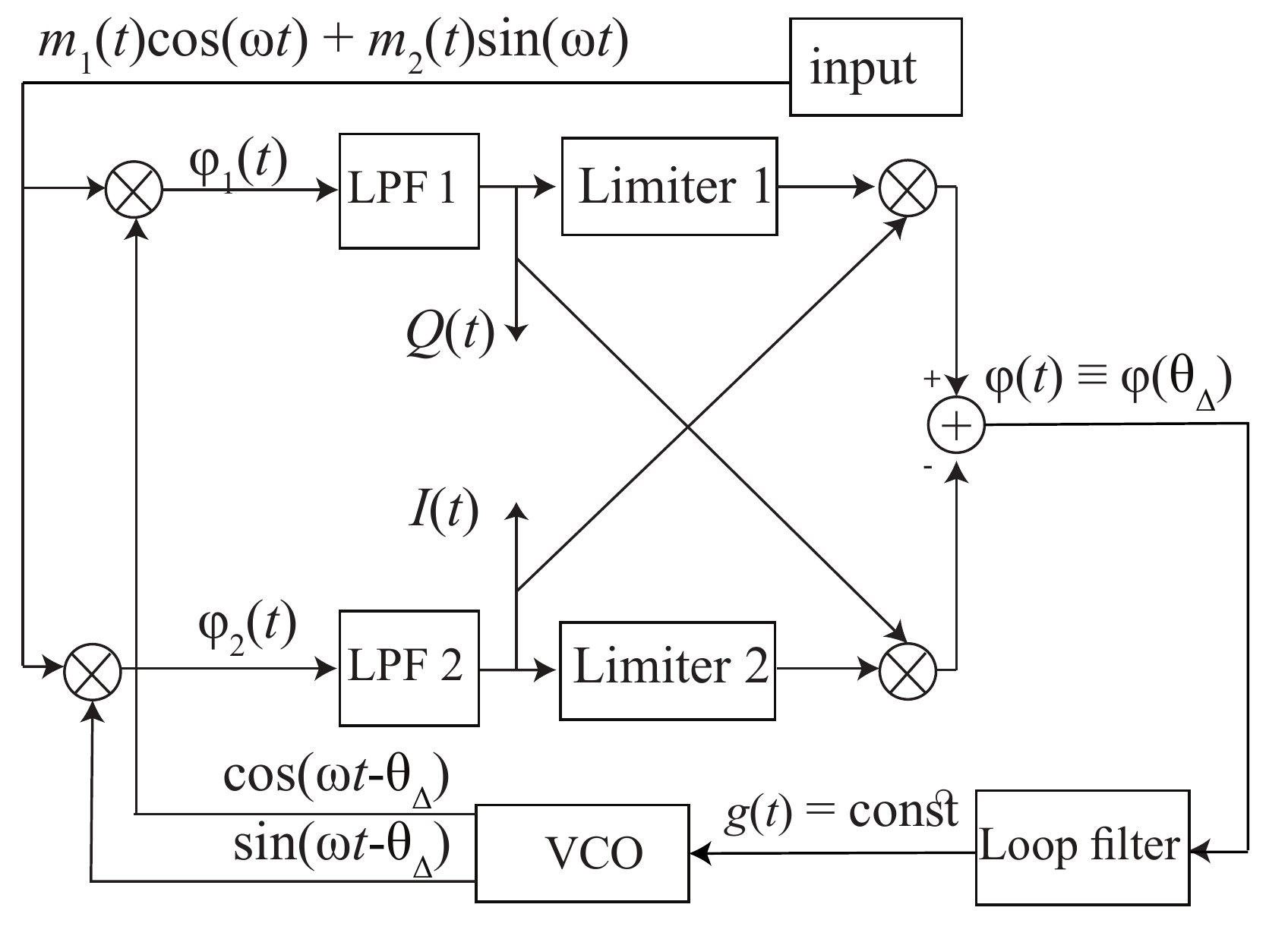}
  \caption{QPSK Costas loop after transient process.}
\label{costas_after_sync}
\end{figure}

Consider the Quadrature Phase Shift Keying Costas loop
(QPSK Costas loop)
after transient processes (see Fig.~\ref{costas_after_sync}).
The input QPSK signal has the form
\begin{equation}
\notag
\begin{aligned}
m_1(t)\cos(\omega t) + m_2(t)\sin(\omega t),
\end{aligned}
\end{equation}
where $m_{1,2}(t) = \pm 1$ is data signal,
$\sin(\omega t)$ and $\cos(\omega t)$ are sinusoidal carriers, $\theta_{\rm ref}(t) = \omega t$ --- phase of input signal.
The VCO has two outputs with $90^{o}$ phase difference: $\cos(\omega t - \theta_{\Delta})$ and $\sin(\omega t - \theta_{\Delta})$, with $\theta_{\rm vco}(t) = \omega t - \theta_{\Delta}$ --- their phase.

After multiplication of VCO signal
and the input signal by multiplier block ($\otimes$) on the upper branch one has
\begin{equation}
  \notag
  \begin{aligned}
      &
      \varphi_1(t) =
      \Big(
        m_1(t)\cos(\omega t)
        + m_2(t)\sin(\omega t)\Big)\cos(\omega t - \theta_{\Delta}).
  \end{aligned}
\end{equation}
On the lower branch  the output signal of VCO is multiplied by the input signal:
\begin{equation}
  \notag
  \begin{aligned}
    &
      \varphi_2(t) = \Big(
        m_1(t)\cos(\omega t)
      + m_2(t)\sin(\omega t)\Big)\sin(\omega t - \theta_{\Delta}).
  \end{aligned}
\end{equation}

{\bf Assumption 1}.
{\it
  The initial states of filters $x_1(0)$, $x_2(0)$, and $x(0)$ do not affect the synchronization of the loop
  (since for the properly designed filters,
  the impact of filter's initial state on its output decays exponentially with time).
}\smallskip

\begin{figure}[H]
  \includegraphics[scale=0.4]{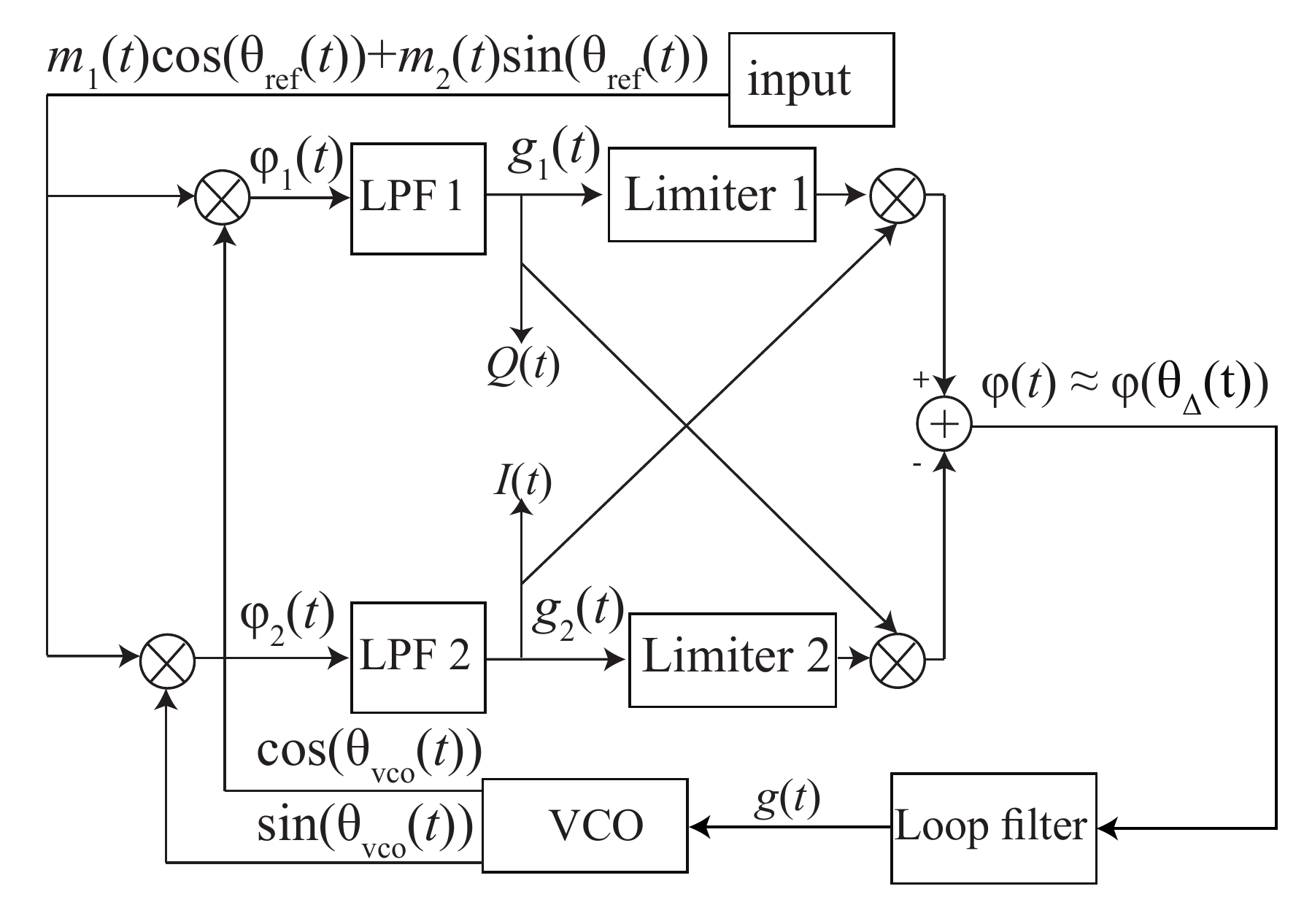}
\caption{QPSK Costas loop before synchronization.}
\label{costas_before}
\end{figure}

Assumption~1 allows one to consider the dependence of the filter output only on its input
ignoring its internal state (see Fig.~\ref{costas_before}).

{\bf Assumption 2}
{\it The terms, whose frequency is about twice the carrier frequency,
     do not affect the synchronization of the loop
     (since they are supposed to be completely suppressed by the low-pass filters).
}
\smallskip

Here, from an engineering point of view,
the high-frequency terms
($\cos(2\omega t - \theta_{\Delta})$ and $\sin(2\omega t - \theta_{\Delta})$
are removed by ideal low-pass filters LPF 1 and LPF 2.
Therefore the consideration of such approximations
doesn't change the outputs of low-pass filter
and is not essential for the analysis of synchronization.

In this case, by Assumption 1
the signals $Q(t)$ and $I(t)$ on the upper and lower branches
can be approximated as
\begin{equation}\label{g1g2-approx}
  \begin{aligned}
    &
      Q(t) \approx
       \frac{1}{2}\Big(m_1(t)\cos(\theta_{\Delta}) + m_2(t)\sin(\theta_{\Delta})
      \Big),
    \\
    &
      I(t) \approx
      \frac{1}{2}\Big( - m_1(t)\sin(\theta_{\Delta}) + m_2(t)\cos(\theta_{\Delta})
      \Big).
  \end{aligned}
\end{equation}
For small values of $\theta_{\Delta}$ we get demodulated data
\begin{equation}\label{g1g2-approx2}
  \begin{aligned}
    &
      Q(t) \approx
  \frac{1}{2}m_1(t),
  \quad
      I(t) \approx
    \frac{1}{2}m_2(t).
  \end{aligned}
\end{equation}

Consider Costas loop before synchronization
(see Fig.~\ref{costas_before})
in the case when the phase difference is not constant:
\begin{equation}\label{thetadelta}
  \theta_{\Delta}(t) = \theta_{\rm ref}(t)-\theta_{\rm vco}(t) \neq const.
\end{equation}

{\bf Caveat to Assumption 2}.
While Assumption~2 is reasonable from a practical point of view,
its use in the analysis of Costas loop
requires further consideration
(see, e.g., \citep{PiqueiraM-2003}).
Here the application of averaging methods
allows one to justify Assumption~2
and obtain the conditions
under which Assumption~2 can be used
(see, e.g., \citep{LeonovKYY-2012-TCASII,LeonovKYY-2016-DAN}).
\smallskip

After the filtration, both signals $\varphi_1(t)$ and $\varphi_2(t)$
pass through the limiters.
Then the outputs of the limiters $\sign\big(Q(t)\big)$ and $\sign\big(I(t)\big)$
are multiplied by $I(t)$ and $Q(t)$.
By Assumption 2 and corresponding formula \eqref{g1g2-approx}
the difference of these signals
\begin{equation}\label{loop-filter-input-approx}
  \begin{aligned}
   &  \varphi(t) = I(t)\sign\big(Q(t)\big) - Q(t)\sign\big(I(t)\big)
  \end{aligned}
\end{equation}
can be approximated as
\begin{equation}\label{phi approx}
  \begin{aligned}
   &
    \varphi(t)
    \approx
    \frac{1}{2}
      \Big(
          - m_1(t)\sin(\theta_{\Delta}(t)) + m_2(t)\cos(\theta_{\Delta}(t))
      \Big)
      \\
      &
      \sign
          \Big(
            m_1(t)\cos(\theta_{\Delta}(t)) + m_2(t)\sin(\theta_{\Delta}(t))
          \Big)
      -
   \\
   &
    -
    \frac{1}{2}
    \Big(
      m_1(t)\cos(\theta_{\Delta}(t)) + m_2(t)\sin(\theta_{\Delta}(t))
    \Big)
    \\
    &
    \sign
      \Big(
          - m_1(t)\sin(\theta_{\Delta}(t)) + m_2(t)\cos(\theta_{\Delta}(t))
      \Big)
    =
   \\
   &
   = \varphi(\theta_{\Delta}(t)) =
   \left\{
      \begin{array}{ll}
        -\sin(\theta_{\Delta}(t)), &  -\frac{\pi}{4}< \theta_{\Delta}(t) < \frac{\pi}{4}, \\
        \cos(\theta_{\Delta}(t)), &  \frac{\pi}{4}< \theta_{\Delta}(t) < \frac{3\pi}{4}, \\
        \sin(\theta_{\Delta}(t)), &  \frac{3\pi}{4}< \theta_{\Delta}(t) < \frac{5\pi}{4}, \\
        -\cos(\theta_{\Delta}(t)), &  \frac{5\pi}{4}< \theta_{\Delta}(t) < -\frac{\pi}{4}. \\
      \end{array}
   \right.
  \end{aligned}
\end{equation}
Here $\varphi(\theta_{\Delta}(t))$ is a piecewise-smooth function.
It should be noted, that function $\varphi(\theta_{\Delta}(t))$
depends on $m_{1,2}$ in points $\theta_{\Delta} = \pm \frac{\pi}{4}, \pm \frac{3\pi}{4}$.
\begin{figure}[H]
  \includegraphics[scale=0.4]{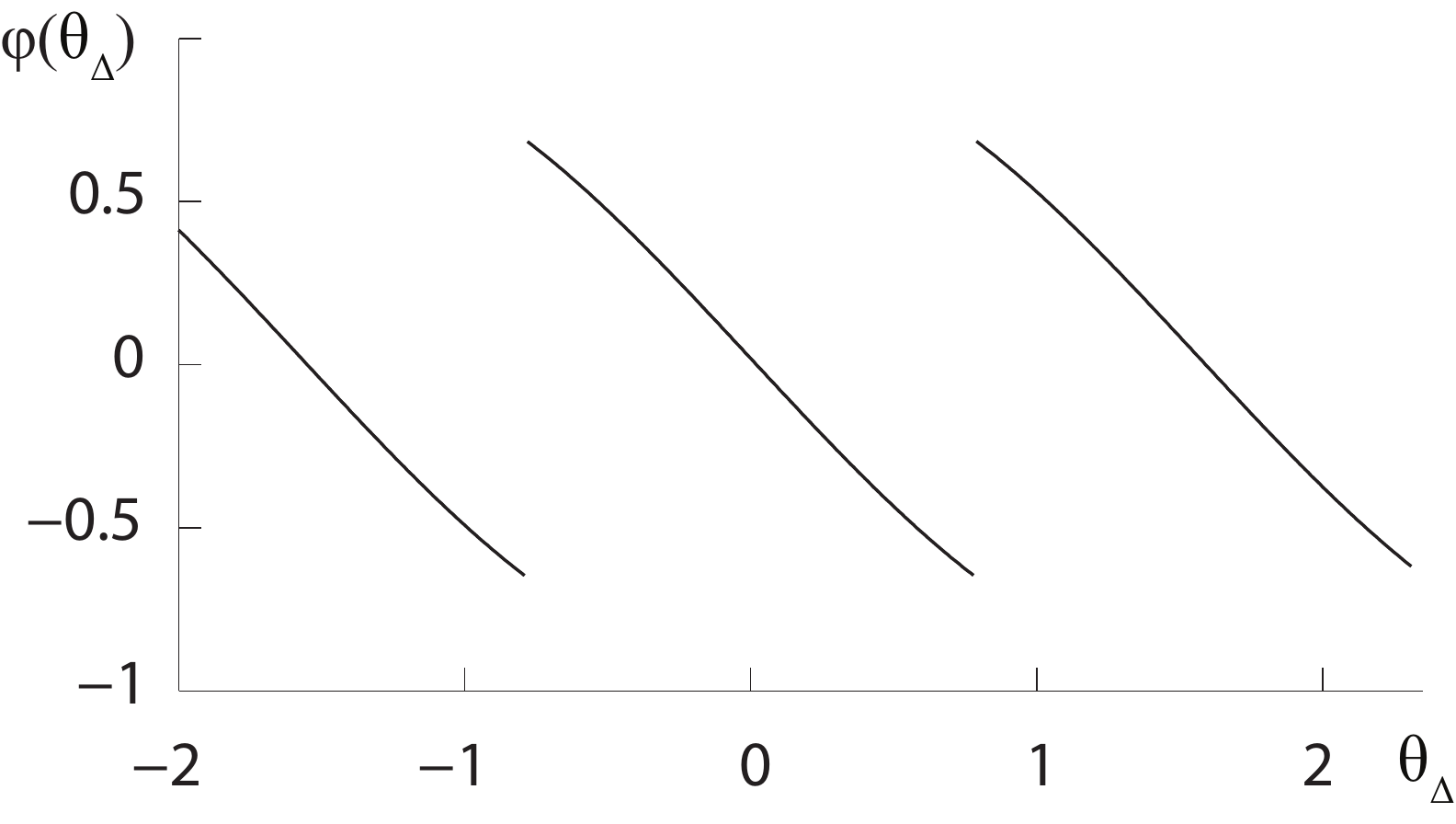}
\caption{Phase detector characteristic of QPSK Costas loop $\varphi(\theta_{\Delta})$.}
\label{costas_out}
\end{figure}
The resulting signal $\varphi(t)$
after the filtration by the loop filter forms the
control signal $g(t)$  for the VCO.

{\bf Assumption 3}.
{\it
The data signals $m_{1,2}(t)$ do not affect the synchronization of the loop.
}

{
Assumptions~1--3 together
lead to the concept of so-called \emph{ideal low-pass filter}.
It removes the upper sideband,
 whose frequency is about twice carrier frequency (Assumption~2),
and passes the lower sideband without change (Assumptions 1,3).
Thus it is assumed that the lower sideband of $\varphi_1(t)$ and $\varphi_2(t)$
are passed without changes and the transmitted data $m_{1,2}(t)$ is neglected
in the signal $\varphi(t)$ (see equation \eqref{phi approx}).
For $m_{1,2}(t) \equiv const$ approximations \eqref{g1g2-approx} depend on the phase difference of signals only, i.e. two multiplier blocks ($\otimes$) on the upper and lower branches
operate as phase detectors.

}

{\bf Caveat to Assumption 3}.
Low-pass filters
can not operate perfectly at moments of changing $m_{1,2}(t)$,
therefore the data pulse shapes are no longer
ideal rectangular pulses after filtration
due to distortion, created by the low-pass filters.
This can lead to incorrect conclusions on the performance of the loop.
One of known examples is so-called false-lock:
while for $m_{1,2}(t) \equiv const$
the loop acquires lock and proper synchronization of the carrier and VCO frequencies,
for time-varying $m(t) \neq const$
the loop can acquire lock without proper synchronization of the frequencies (false lock)
\citep{Olson-1975,Simon-1978,Lindsey-1978}.
To avoid such undesirable situation one may
try to choose loop parameters in such a way that
the synchronization time is less than the time between changes in the data signal $m_{1,2}(t)$
or to modify the loop design (see, e.g., \citep{Olson-1975}).
Another way is to perform the nonlinear nonlocal analysis of the loop
(see, e.g., \citep{Stensby-1989,Stensby-2002})
to identify unsuitable parameters.
\smallskip

{\bf Caveat to Assumption 1}.
If in Fig.~\ref{costas_before} the loop is out of lock,
i.e. synchronization is not achieved,
filters' initial states cannot be ignored
and must be taken into account.
Really, low-pass filters with nonzero initial states may change
the lower sideband (see expressions \eqref{g1g2-approx})
and affect the synchronization of the loop.
For rigorous consideration of low-pass filters
one has to use mathematical models of filters
instead of approximations \eqref{g1g2-approx}.
Since the low-pass filters LPF 1 and LPF 2 are mostly used
for data demodulation, the effect of nonzero initial state of filter
on transient processes will be discussed for the loop filter,
which is used to provide synchronization.
\smallskip

The relation between the input $\varphi(t)$
and the output $g(t)$ of the Loop filter has the form
\begin{equation}\label{loop-filter}
 \begin{aligned}
 & \frac{dx}{dt} = A x + b \varphi(t),
 \ g(t) = c^*x + h\varphi(t).
 \end{aligned}
\end{equation}
Here $A$ is a constant matrix,
vector $x(t)$ is a filter state,
$b,c$ are constant vectors,
and $x(0)$ is initial state of filter.
The solution of equation \eqref{loop-filter} with initial data $x(0)$
(filter initial state) is as follows
\begin{equation}\label{loop-filter-solution}
 \begin{array}{c}
 g(t) = \alpha_0(t) + h \varphi(t)+
 \int\limits_0^t
 \gamma(t - \tau)\varphi(\tau)
 {\rm d}\tau.
 \end{array}
\end{equation}
Here $\gamma(t - \tau)=c^*e^{A(t-\tau)}b$ is an impulse response function of filter
and $\alpha_0(t) = \alpha_0(t,x(0))= c^*e^{At}x(0)$ is an exponentially damped function
(i.e. the matrix A is stable).
Corresponding transfer function takes the form\footnote{
In the control theory \citep{LeonovK-2014-book}
it is defined with opposite sign:  $c^*(A-sI)^{-1}b-h$}
\begin{equation}
  H(s) = -c^{*}(A - sI)^{-1}b + h.
\end{equation}

The control signal $g(t)$ is used to adjust VCO frequency to the frequency of input carrier signal
\begin{equation} \label{vco first}
   \dot\theta_{\rm vco}(t) = \omega_{\rm vco}(t) = \omega_{\rm vco}^{\text{free}} + K_{\rm vco}g(t).
\end{equation}
Here $\omega_{\rm vco}^{\text{free}}$ is free-running frequency of VCO
and $K_{\rm vco}$ is VCO gain.
Note that the initial VCO frequency (at $t=0$) is as follows
\begin{equation}
  \omega_{\rm vco}(0) = \omega_{\rm vco}^{\text{free}} + K_{\rm vco}\alpha_0(0)+K_{\rm vco}h\varphi(\theta_{\Delta}(0))
  \neq \omega_{\rm vco}^{\text{free}}.
\end{equation}

If the frequency of input carrier is a constant
\begin{equation}\label{omega1-const}
   \dot\theta_{\rm ref}(t) = \omega_{\rm ref}(t) \equiv \omega_{\rm ref},
\end{equation}
then equations \eqref{loop-filter-input-approx}-\eqref{vco first}
give the following mathematical model of Costas loop
\begin{equation} \label{mathmodel-class}
 \begin{aligned}
   & \dot\theta_{\Delta} =
   \omega_{\rm ref} - \omega_{\rm vco}^{\rm free}-K_{\rm vco}\alpha_0(t)- \\
   & -K_{\rm vco}
   \Big(
    h\varphi(\theta_{\Delta}) +\int\limits_0^t
    \gamma(t - \tau)\varphi(\theta_{\Delta}(\tau)){\rm d}\tau
   \Big). \\
 \end{aligned}
\end{equation}

{\bf Assumption 4 (Corollary 1 of Assumption 1)}.
{\it
  Free output of loop filter $\alpha_0(t)$ does not affect the synchronization of the loop
  since $\alpha_0(t)$ is an exponentially damped function.
}\smallskip

\begin{figure}[H]
  \includegraphics[scale=0.6]{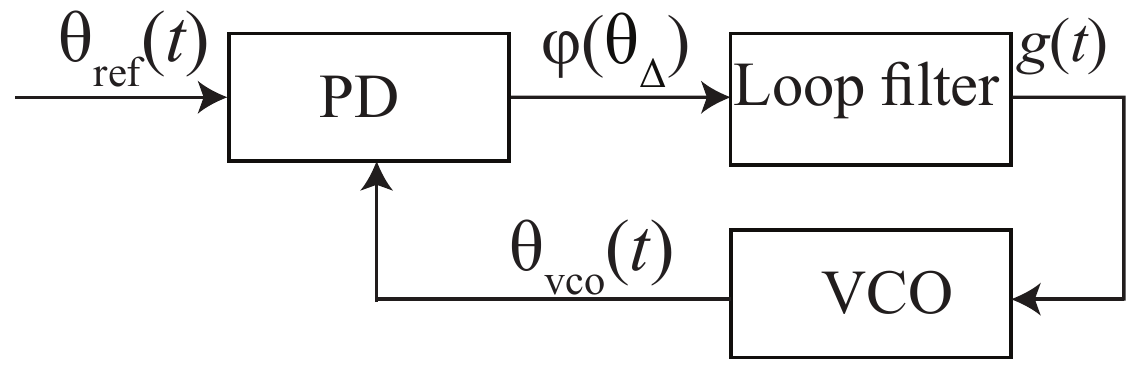}
  \caption{Classical mathematical model of QPSK Costas loop.}
  \label{pll-qpsk}
\end{figure}

For $h=0$ Assumption~4 allows one to obtain
the classical mathematical model of Costas loop (see Fig.~\ref{pll-qpsk})
\begin{equation} \label{mathmodel-class-simple}
 \begin{aligned}
   & \dot\theta_{\Delta} =
   \omega_{\Delta}^{\rm free}
   - K_{\rm vco}\int\limits_0^t
   \gamma(t - \tau)\varphi(\theta_{\Delta}(\tau)){\rm d}\tau \\
   & \theta_{\Delta}(t) = \theta_{\rm ref}(t)-\theta_{\rm vco}(t), \
   \omega_{\Delta}^{\rm free} = \omega_{\rm ref}-\omega_{\rm vco}^{\text{free}}
 \end{aligned}
\end{equation}

{\bf Caveat to Assumption 4}.
For high-order filter,
two different initial states $\tilde{x}(0)$ and $\tilde{\tilde{x}}(0)$ may lead
to identical values of $\alpha_0(0,\tilde{x}(0))=\alpha_0(0,\tilde{\tilde{x}}(0))$
but different functions $\alpha_0(t,\tilde{x}(0))$ and $\alpha_0(t,\tilde{\tilde{x}}(0))$
(to avoid this effect it is necessary to assume the observability of system \eqref{loop-filter}).

Since nonlinear mathematical model of Costas loop \eqref{mathmodel-class-simple}
is hard to analyze, in practice, for its analysis
it is widely used numerical simulation and linearization.
In the case when the phase difference of signals is small
one can consider a linearized mathematical model of Costas loop,
using the linearization $\varphi(\theta_{\Delta}) \approx K\theta_{\Delta}$.
This allows one to estimate hold-in range by the same
methods that were developed for analysis and design of classical PLLs
(see, e.g., \citep{Gardner-1966,Viterbi-1966,Lindsey-1972,Shahgildyan-1972},
and others).
Linearized model \eqref{mathmodel-class},
where $\varphi(\Delta\theta)$ is changed by $K\theta_{\Delta}$,
may be used for analysis in the case when the loop is in lock,
but  analysis of the acquisition behavior cannot be accomplished
using linearized models.

Next we discuss rigorous derivation of nonlinear mathematical model.
The relation between the inputs $\varphi_{1,2}(t)$
and the outputs $g_1(t)=Q(t)$ and $g_2(t)=I(t)$
of the low-pass filters is similar to \eqref{loop-filter}:
\begin{equation}\label{LPF}
 \begin{aligned}
 & \frac{dx_{1,2}}{dt} = A_{1,2} x_{1,2} + b_{1,2}\varphi_{1,2}(t),
 \ g_{1,2}(t) = {c_{1,2}}^*x_{1,2}.
 \end{aligned}
\end{equation}
Here $A_{1,2}$ are constant matrices,
the vectors $x_{1,2}(t)$ are filter states,
$b_{1,2},c_{1,2}$ are constant vectors,
and $x_{1,2}(0)$ are initial states of filters.

\begin{figure}[H]
\centering
  \includegraphics[scale=0.4]{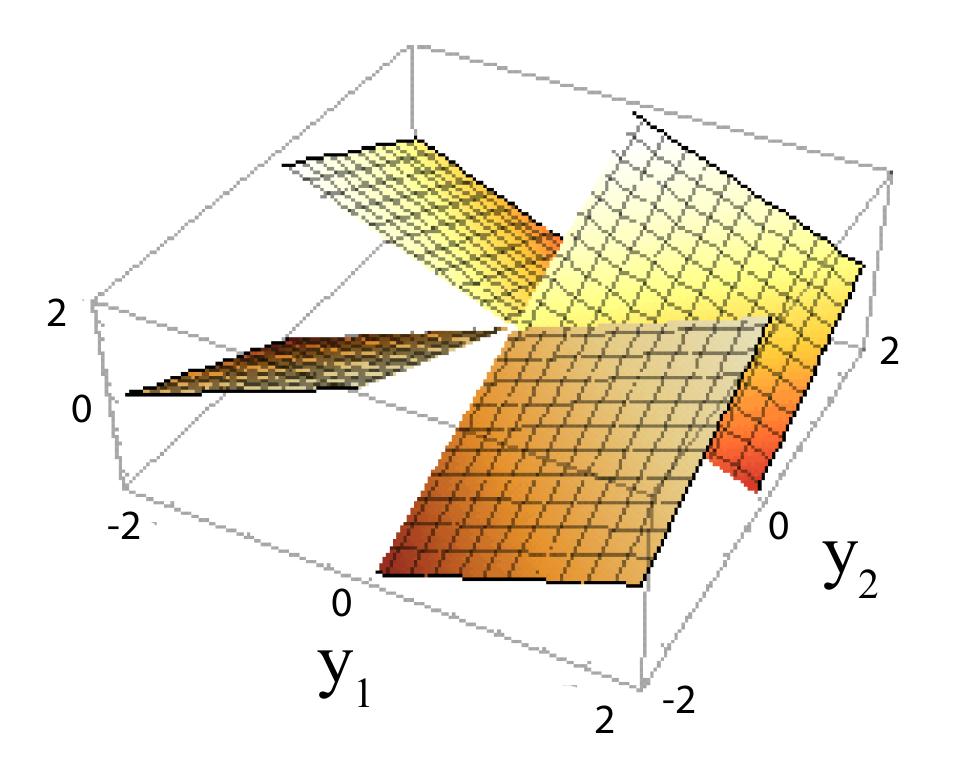}
\caption{Plot of the function $y_1\sign(y_2)- y_2\sign(y_1)$.}
\label{function plot}
\end{figure}

Then, taking into account \eqref{LPF}, \eqref{loop-filter}, and \eqref{vco first},
one obtains \emph{mathematical model in the signal space}
describing \emph{physical model} of QPSK Costas loop:
\begin{equation} \label{prev diff eq}
    \begin{aligned}
        & \dot{x_1} = A_1 x_1 + \\ & + b_1\cos(\theta_{\rm vco})
          \big(m_1(t)\cos(\theta_{\rm ref}(t)) + m_2(t)\sin(\theta_{\rm ref}(t))\big), \\
        & \dot{x_2} = A_2 x_2 + \\ & + b_2\sin(\theta_{\rm vco})
          \big(m_1(t)\cos(\theta_{\rm ref}(t)) + m_2(t)\sin(\theta_{\rm ref}(t))\big), \\
        & \dot{x} = A x + b(\sign(c_1^*x_1)(c_2^*x_2)-\sign(c_2^*x_2)(c_1^*x_1)), \\
        & \dot\theta_{\rm vco} = \omega_{\rm vco}^{\text{free}} + K_{\rm vco}(c^*x) + \\
        & \qquad + K_{\rm vco}\big((c_2^*x_2)\sign(c_1^*x_1)-(c_1^*x_1)\sign(c_2^*x_2)\big).
    \end{aligned}
\end{equation}
Here $\theta_{\rm vco}(0)$ is the initial phase shift of VCO
and the vectors $x_{1,2}(0), x(0)$ are initial states of filters
(so Assumptions 2 and 4 are not used).
Thus the initial VCO frequency (at $t=0$) has the form
\begin{equation} \label{freq-init-ss}
  \begin{aligned}
   &\omega_{\rm vco}(0)=\omega_{\rm vco}^{\text{free}} + K_{\rm vco}c^*x(0) + \\
   &
   + K_{\rm vco}\big((c_2^*x_2(0))\sign(c_1^*x_1(0))-(c_1^*x_1(0))\sign(c_2^*x_2(0))\big).
    \end{aligned}
\end{equation}
Right-hand side of \eqref{freq-init-ss} is discontinuous.
Fig.~\ref{function plot} shows function $y_1\sign(y_2)- y_2\sign(y_1)$ around zero.

The mathematical model in signal space \eqref{prev diff eq}\
is nonlinear nonautonomous discontinuous differential system,
so in general case its analytical study is a difficult task
even for the continuous case when $m_{1,2}(t)\equiv const$.
Moreover it is a slow-fast system, so its numerical study
is rather complicated for the high-frequency signals.
The problem is that it is necessary
to consider simultaneously both very fast time scale
of the signals $sin(\theta_{1,2}(t))$
and slow time scale of phase difference between the signals
$\theta_{\Delta}(t)$,
therefore one very small simulation time-step
must be taken over a very long total simulation period
\citep{Goyal-2006,KuznetsovKLNYY-2014-ICUMT-QPSK}.

To overcome these problems in PLL and classic Costas loop,
in place of using Assumption~2
one can apply averaging methods \citep{KrylovB-1947,MitropolskyB-1961,Samoilenko-2004-averiging,SandersVM-2007-averaging}
and consider \emph{a simplified mathematical model in the signal's phase space}.
Remark that  classical averaging approach requires Lipschitz condition,
which is not satisfied in the case of QPSK Costas loop system.


It is useful to formalize engineering Assumption~2 and
the explanation of low-pass filters operation
in the following way
\begin{equation}
\label{lpf_conditions}
\begin{aligned}
 &
 \int_{t_0}^t \gamma_{1,2}(t-\tau)\sin\theta(\tau)d\tau
 = \\
 & =
 \sin\big(\theta(t)\big) + O(\frac{1}{\omega_{min}}),
 \ \forall \dot \theta(t) < \omega_{\Delta}^{\text{max}},
 \\
  &
 \int_{t_0}^t \gamma_{1,2}(t-\tau)\sin\theta(\tau)d\tau
 = O(\frac{1}{\omega_{min}}),
 \\
 &
  \forall \dot \theta(t) > \frac{C}{\sqrt \omega_{min}},
  \\
 & \gamma_{1,2}(t - \tau) = c_{1,2}^*e^{A_{1,2}(t - \tau)}b_{1,2}, \\
 & \alpha_{1,2}(t) = c_{1,2}^*e^{A_{1,2} t}x_{1,2}(0).
\end{aligned}
\end{equation}
Here
$\omega_{\text{vco, ref}}(t)>\omega^{\rm min}>0$,
$1/\sqrt{\omega^{\rm min}} < |\omega_{\text{ref}}(t)-\omega_{\text{vco}}(t)|
<\omega_\Delta^{\rm max}$ on sufficiently large time interval,
$t_0$ is a moment of time such that $\alpha_{1,2}(t) \leq \frac{1}{\omega^{\rm min}}$ for $t > t_{0}$.
If the initial states of low-pass filters LPF 1 and LPF 2 are zero,
then $\alpha_1(t) = \alpha_2(t) = 0$ (see Assumption~1, Assumption~4).

Applying \eqref{lpf_conditions} to \eqref{prev diff eq}, one obtains
\begin{equation}
  \label{averaged system with O}
  \begin{aligned}
    & \dot{x} = A x + b \varphi(\theta_{\Delta}),
    \\
    & \dot\theta_{\Delta} = \omega_{\Delta}^{\rm free} - K_{\rm vco}(c^*x) - K_{\rm vco}h\varphi(\theta_{\Delta}),
    \\
    & \varphi(\theta_{\Delta}) =
    \frac{1}{\sqrt{2}}\bigg(
        \sin(\theta_{\Delta}(\tau)+\frac{\pi}{4})
        \\
        & \qquad\qquad
        \sign\Big(
              \sin(\theta_{\Delta}(t)-\frac{\pi}{4})
             \Big)
        -
        \\
        & \qquad\qquad
        -
        \sin(\theta_{\Delta}(\tau)-\frac{\pi}{4})
        \\
        & \qquad\qquad
        \sign\Big(
                \sin(\theta_{\Delta}(t)+\frac{\pi}{4})
              \Big)
    \bigg) + O({1 \over \sqrt\omega_{min}}).
  \end{aligned}
\end{equation}
Corresponding detailed discussion can be found in \citep{LeonovKYY-2016-DAN}.

{\bf Assumption 5 (Corollary of Assumptions 1-3)}.
{\it Solutions of system \eqref{prev diff eq}
under condition \eqref{lpf_conditions}
are close to the solutions of the following system
(i.e. $O(\frac{1}{\sqrt{\omega_{min}}})$ can be neglected)}
\begin{equation}
 \label{final_system}
 \begin{aligned}
 & \dot{x} = A x + b \varphi(\theta_{\Delta}), \\
 & \dot\theta_{\Delta} = \omega_{\Delta}^{\rm free} - K_{\rm vco}c^*x - K_{\rm vco}h\varphi(\theta_{\Delta}), \\
    & \varphi(\theta_{\Delta}) =
    \frac{1}{\sqrt{2}}\bigg( \\ &
          \sin(\theta_{\Delta}(\tau)+\frac{\pi}{4})
          \sign\Big(
                  \sin(\theta_{\Delta}(t)-\frac{\pi}{4})
                \Big)
          -
          \\
          &
          -
          \sin(\theta_{\Delta}(\tau)-\frac{\pi}{4})
          \sign\Big(
                  \sin(\theta_{\Delta}(t)+\frac{\pi}{4})
                \Big)
        \bigg).
 \end{aligned}
\end{equation}
Here function $\varphi(\theta_{\Delta})$ is a
phase detector characteristic of QPSK Costas loop for sinusoidal signals,
which is used in classical books.
Note that here the phase detector operation include
operations of multipliers, limiters, LPF 1, and LPF 2.

Let us determine equilibrium points.
Consider a transfer function of the capacitor-based filter without parasitic resistance
 \begin{equation}
 \label{Cs}
 \begin{aligned}
 & F(s) = \frac{1}{Cs}.
 \\
 \end{aligned}
 \end{equation}
 System \eqref{final_system} with this filter takes the following form
 \begin{equation}
 \begin{aligned}
 & \dot{x} = {1 \over C} \varphi(\theta_{\Delta}), \\
 & \dot\theta_{\Delta} = \omega_{\Delta}^{\rm free} - Lx.
 \end{aligned}
 \end{equation}
 This system has the following equilibrium points
 \begin{equation}
 \begin{aligned}
 & x = {\omega_{\Delta}^{\rm free} \over L},
 \quad \varphi(\theta_{\Delta}) = 0.\\
 \end{aligned}
 \end{equation}

If $A$ is a non-singular matrix, then equilibrium points are determined by
the following system
\begin{equation}
 \begin{aligned}
 & x = A^{-1} b \varphi(\theta_{\Delta}), \\
 & \varphi(\theta_{\Delta}) = {\omega_{\Delta}^{\rm free} \over (K_{\rm vco}c^*A^{-1}b + K_{\rm vco}h)}. \\
 \end{aligned}
\end{equation}
Denote
\begin{equation}
\begin{aligned}
& \gamma = {\omega_{\Delta}^{\rm free} \over (K_{\rm vco}c^*A^{-1}b + K_{\rm vco}h)}.
\end{aligned}
\end{equation}
Then
\begin{figure}[H]
  \includegraphics[scale=0.4]{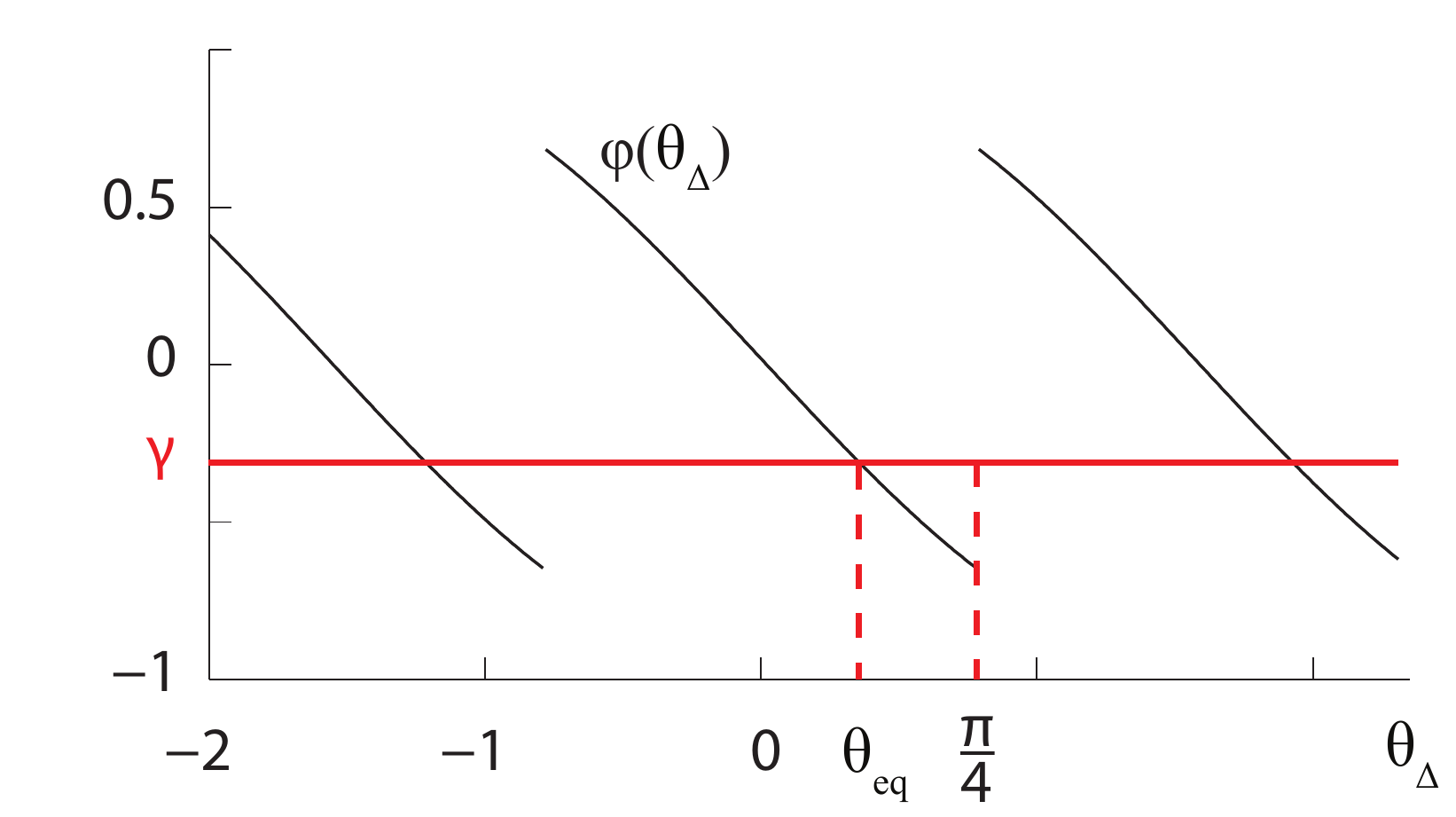}
\caption{Equilibrium points of QPSK Costas loop with filter without $0$ poles.}
\label{qpsk-equilibrium}
\end{figure}
\begin{equation}
\begin{aligned}
& \theta_{eq} = -\arcsin({\omega_{\Delta}^{\rm free} \over L(c^*A^{-1}b + h)}) + \pi k.
\end{aligned}
\end{equation}
The linearized system in the neighborhood of $\theta_{eq}$ is as follows
\begin{equation}
\begin{aligned}
& \dot x = Ax - b\cos(\theta_{eq})\theta_{\Delta - eq},\\
& \dot \theta_{\Delta - eq} = \omega_{\Delta}^{\rm free} - K_{\rm vco}c^*x + K_{\rm vco}\cos(\theta_{eq})\theta_{\Delta - eq}, \\
& \theta_{\Delta - eq} = \theta_{\Delta} - \theta_{eq}
\end{aligned}
\end{equation}
Using equality
\begin{equation}
\begin{aligned}
& \det\left(\begin{array}{cc}
    A & B \\
    C & D \\
  \end{array}\right) = \det A \cdot \det(D - CA^{-1}B),
\end{aligned}
\end{equation}
one can obtain the characteristic polinomial
\begin{equation}
\begin{aligned}
&
  \chi(s) =
  \det\left(\begin{array}{cc}
    A - sI& -b\cos(\theta_{eq}) \\
    -K_{\rm vco}c^* & K_{\rm vco}\cos(\theta_{eq}) - s \\
  \end{array}\right)
\\
&
  = \det(A - sI)
\\
&
    (K_{\rm vco}\cos(\theta_{eq}) - s
       - K_{\rm vco}c^* (A - sI)^{-1} b\cos(\theta_{eq}))
  =
\\
&
  = \det(A - sI)
    (K_{\rm vco}\cos(\theta_{eq}) - s
       + LH(s)\cos(\theta_{eq})).
\end{aligned}
\end{equation}
Denote a filter transfer function $H(s) = {M(s) \over N(s)}$.
Then
\begin{equation}
\begin{aligned}
&
  \chi(s) = -N(s)(K_{\rm vco}\cos(\theta_{eq}) - s) + M(s)L\cos(\theta_{eq})
\end{aligned}
\end{equation}
For the properly designed QPSK Costas Loop $\theta_{eq}$ is a stable point.
If all of the zeros of the characteristic polynomial $\chi(s)$
have negative real parts, then $\theta_{eq}$ is asymptotically stable
equilibrium point.
Thus, the parameters of the filters ($N$, $M$, $L$, and $h$)
 should satisfy this rule.

Suppose that
\begin{equation}
\label{eqb}
  -\frac{\pi}{4} + \pi n < \theta_{\Delta}(t) < \frac{\pi}{4} + \pi n.
\end{equation}
Then the function $\varphi(t)$ becomes continuous, so it is possible to apply
classical averaging theorem. 
Until the initial frequency difference is sufficiently large
and \eqref{eqb} is not satisfied, one has to apply different approach
to investigation of transient processes \citep{LeonovKYY-2016-DAN}.

{\bf Caveat to Assumption 5}.
For rigorous justification of Assumption~5
one has to prove
that $O(\frac{1}{\sqrt{\omega_{min}}})$
does not affect the bahaviour of Limiters.

\section{Counterexamples to the Assumptions}
Note once more that various simplifications
and the analysis of linearized models of control systems
may result in incorrect conclusions\footnote{
see also counterexamples to the filter hypothesis,
Aizerman's and Kalman's conjectures on the absolute stability of nonlinear control systems \citep{KuznetsovLS-2011-IFAC,BraginVKL-2011,LeonovK-2013-IJBC},
and the Perron effects of the largest Lyapunov exponent sign inversions \citep{KuznetsovL-2005,LeonovK-2007}, etc.}.
At the same time the application of nonlinear methods for the analysis of PLL-based models
are quite rare (see, e.g.,
\citep{Abramovitch-1990,ChangTC-1993,Stensby-1997,ShirahamaFYT-1998,Watada-1998,Hinz-2000,Wu-2002,PiqueiraM-2003,SuarezQ-2003,Margaris-2004,VendelinPR-2005-book,banerjee2006,KudrewiczW-2007,Wang-2008,Piqueira-2010-PLL-network,Wiegand-2010,Stensby-2011,Suarez-2012,SarkarSB-2014,ChiconeH-2013,YoshimuraIMM-2013,BestKLYY-2014-DCNPS}).
Further examples demonstrate that the use of Assumptions 1-4
requires further study and rigorous justification.
The following examples are shown
that for the same parameters
the operations of real \emph{physical model} of Costas loop and
\emph{mathematical or physical simplified model},
taking into account one of the above Assumptions,
may differ considerably.

%

\begin{figure}[H]
  \includegraphics[scale=0.5]{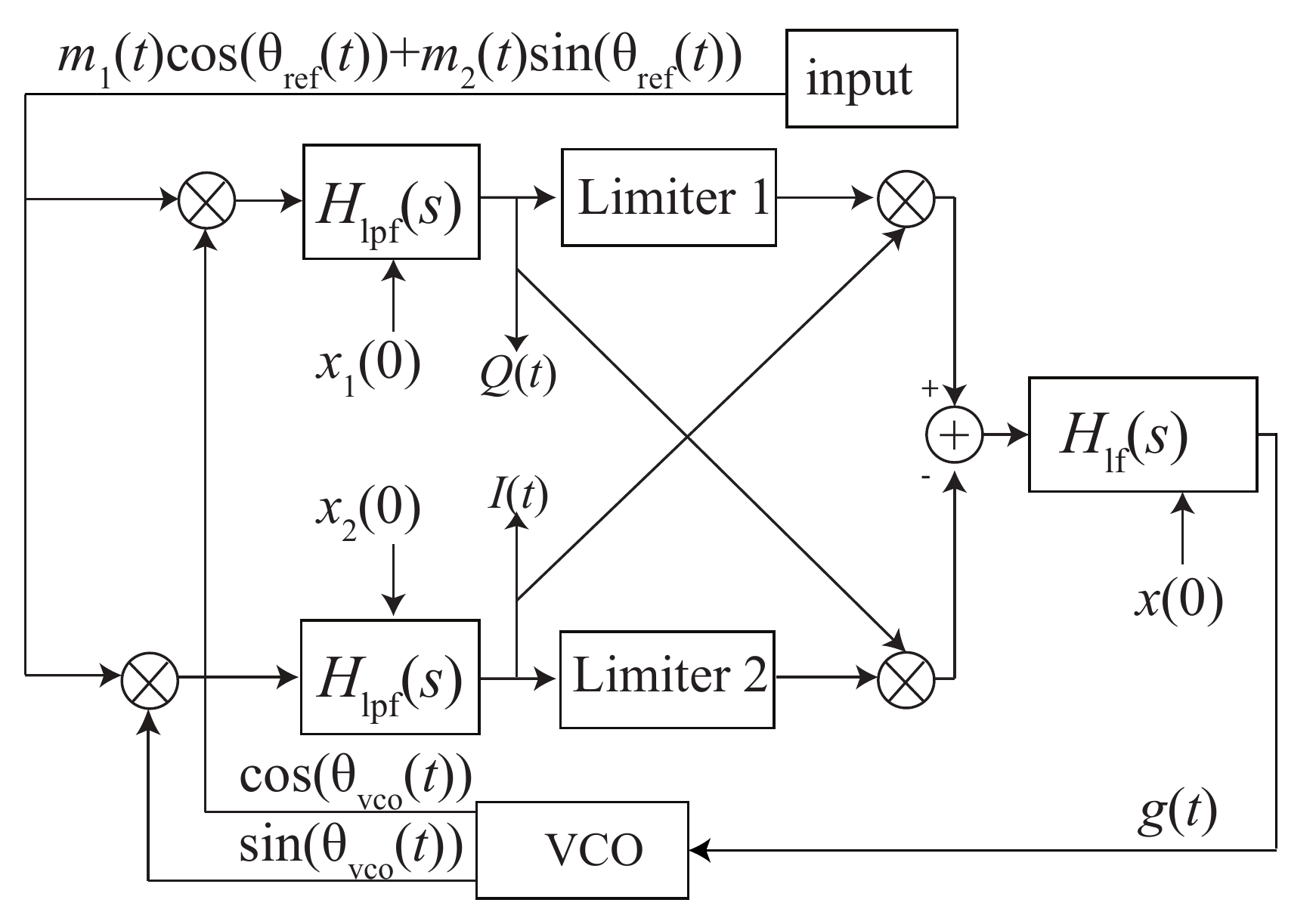}
\caption{QPSK Costas loop width filters defined by their transfer functions.}
\label{qpsk transfer}
\end{figure}

{\bf Simulation model and parameters.}
In engineering practice one of the most popular way to describe linear filter
is considering its transfer function $H(s)$ (see Fig.~\ref{qpsk transfer}).
In the following examples we use
loop filter transfer functions $H_{lf}(s) = \frac{\tau_2 s + 1}{\tau_1 s}$, $$ \tau_1 = 20\cdot10^{-6}, \tau_2 = 4\cdot 10^{-6},$$ described by the equations
  \begin{equation}
    \begin{aligned}
      & \dot x = \xi,\\
      & \sigma = \frac{1}{\tau_1}x + \frac{\tau_2}{\tau_1}\xi,\\
    \end{aligned}
  \end{equation} where $\xi(t)$ is an input of the filter and $\sigma(t)$ is an output of the filter.
  Low pass filters transfer function is $H_{lpf}(s) = \frac{1}{s/\omega_{lpf}+1}$, $\omega_{lpf} = 1.2566\cdot10^6$, and the corresponding equations are
  \begin{equation}
    \begin{aligned}
      & \dot x_{1,2} = -\omega_{lpf} x_{1,2} + \xi,\\
      & \sigma = x_{1,2},\\
    \end{aligned}
  \end{equation}
where carrier frequency is $\omega_{\rm ref}=2\cdot\pi\cdot400000$,
VCO input gain is $6.3165\cdot10^5$;
VCO phase shift is zero;
$m_2(t) = 1$;

\begin{example}
In Fig.~\ref{example_ics_lpf}
is shown that Assumptions~1 and 4 may not be valid:
while physical model with zero initial states of low-pass filters
acquire lock (black color),
physical model with nonzero initial states of low-pass filters
is out of lock (red color).
It should be noted that in Fig.~\ref{example_ics_lpf} initial frequencies of VCO corresponding to the red curve and black curve are the same.
In Fig.~\ref{example_ics_lf} similar example is presented for
nonzero initial state of loop filter.

\begin{figure}[h]
  \includegraphics[width=0.4\textwidth]{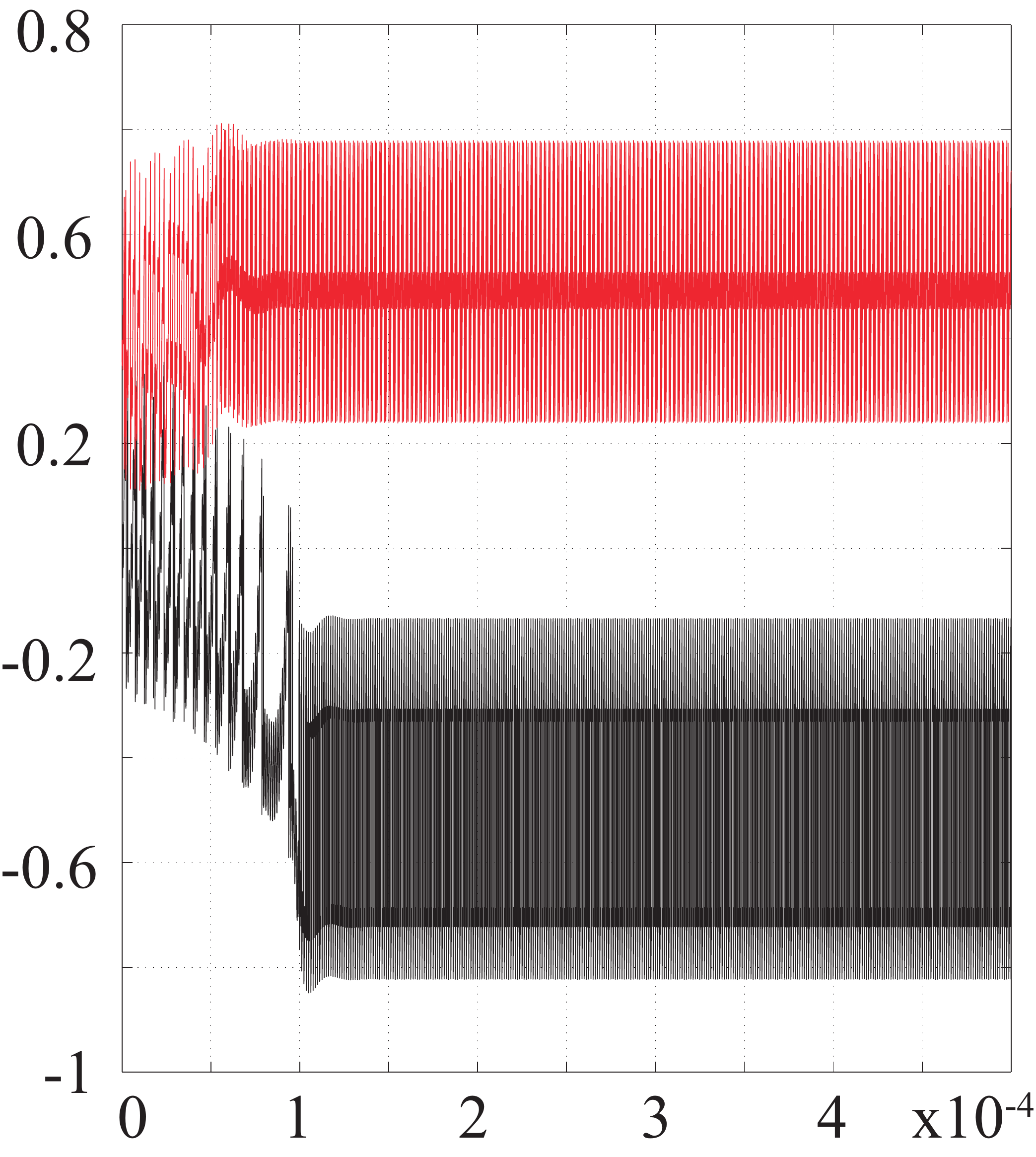}
  \caption{Loop filter output: $m_1(t) = 1$;
  VCO free-running frequency is $2.6314\cdot10^6$ rad/s;
  VCO phase shift is zero;
  $\omega_{lpf} = 1.2566\cdot10^6$;
  initial loop filter state is $0.4$ (red curve) and zero (black curve)}
  \label{example_ics_lf}
\end{figure}

\begin{figure}[h]
  \includegraphics[width=0.4\textwidth]{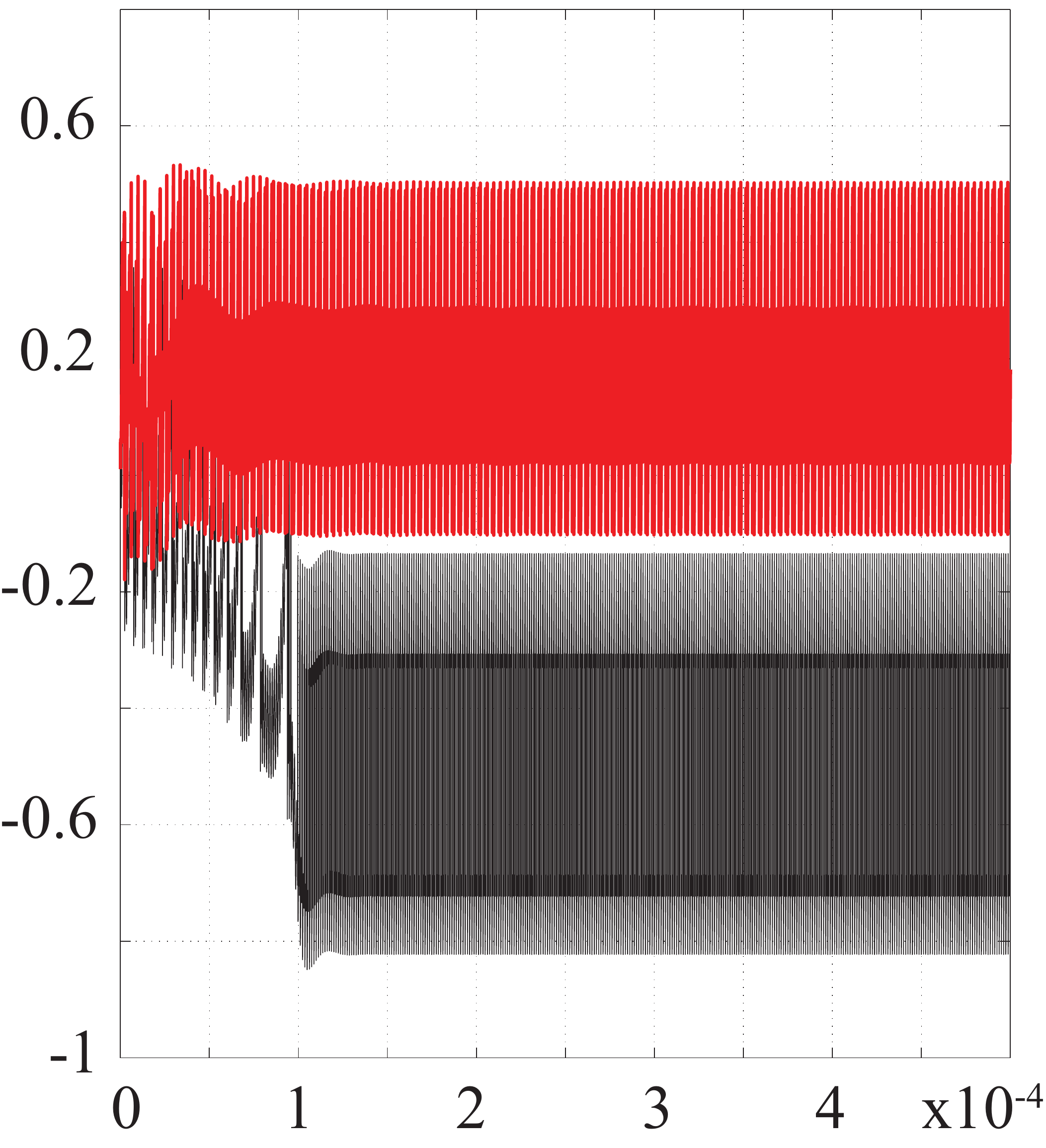}
  \caption{Loop filter output: $m_1(t) = 1$;
  VCO free-running frequency $2.8283\cdot10^6$ rad/s;
  initial conditions of low-pass filters are $30$ for red curve and zero for black curve;
  initial condition for loop filter is zero;
  $\omega_{lpf} = 1.2566\cdot10^6$;
  VCO phase shift is zero}
  \label{example_ics_lpf}
\end{figure}
\end{example}

\begin{example}
In Fig.~\ref{example_pd_instability}
is shown that Assumption~2 may not be valid:
while averaged model 
acquire lock (black), physical model is out of lock (red).

\begin{figure}
  \includegraphics[width=0.4\textwidth]{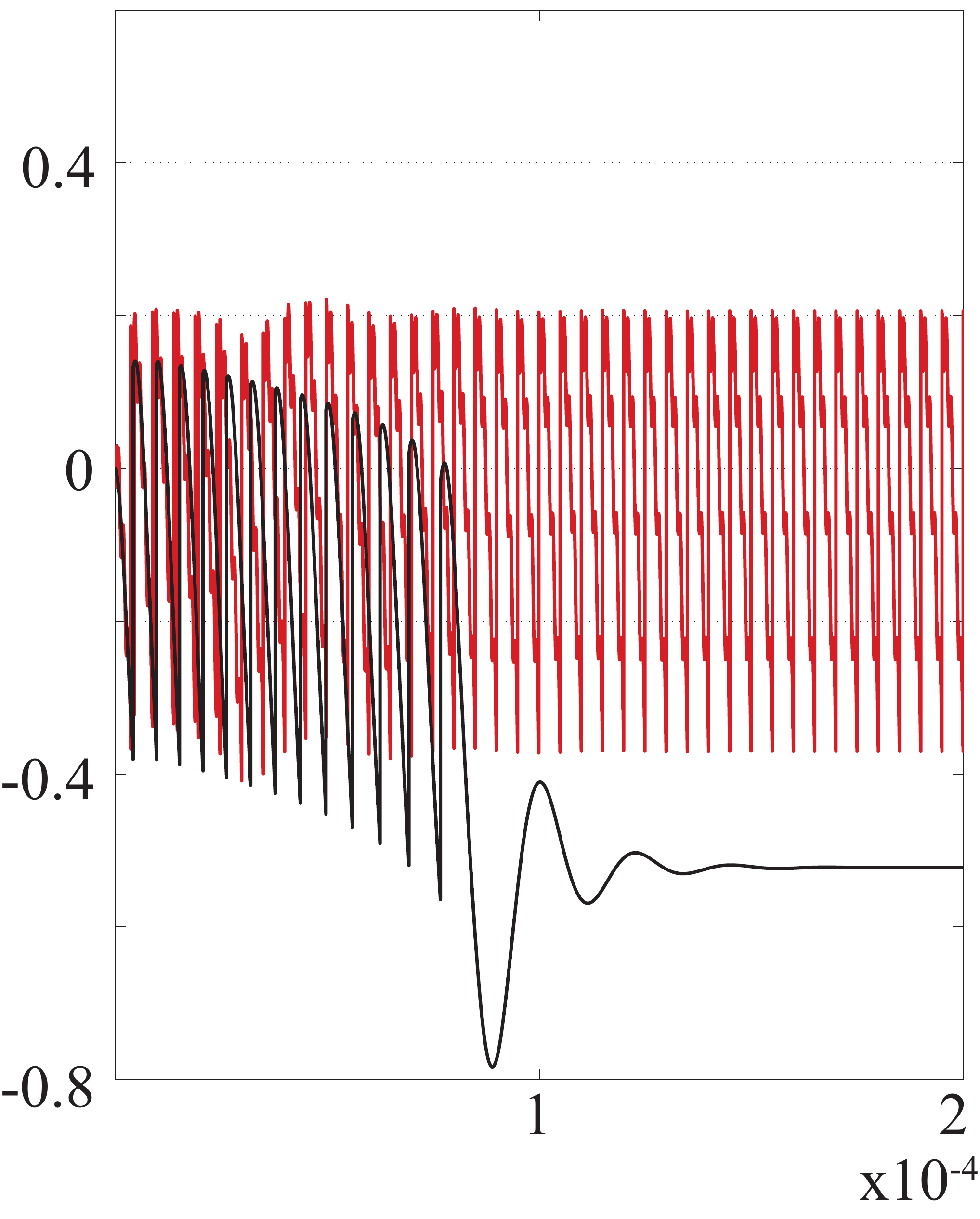}
  \caption{Loop filter output:
  $m_1(t) \equiv 1$;
  Initial conditions of loop filter are zero;
  VCO phase shift is zero;
  VCO free-running frequency is $2.8433\cdot10^6$ rad/s;
  $\omega_{lpf} = 6.2832\cdot10^5$;
  red curve - taking into account signals with twice carrier frequency, black curve - taking into account only low frequency signals.
  }
  \label{example_pd_instability}
\end{figure}
\end{example}

\begin{example}
In Fig.~\ref{example_lock_unlock_with_data}
is shown that Assumption~3 may not be valid:
while physical model (black)
with constant data signal acquire lock,
physical model (red) with nonconstant data signal is out of lock.

\begin{figure}
  \centering
  \includegraphics[width=0.4\textwidth]{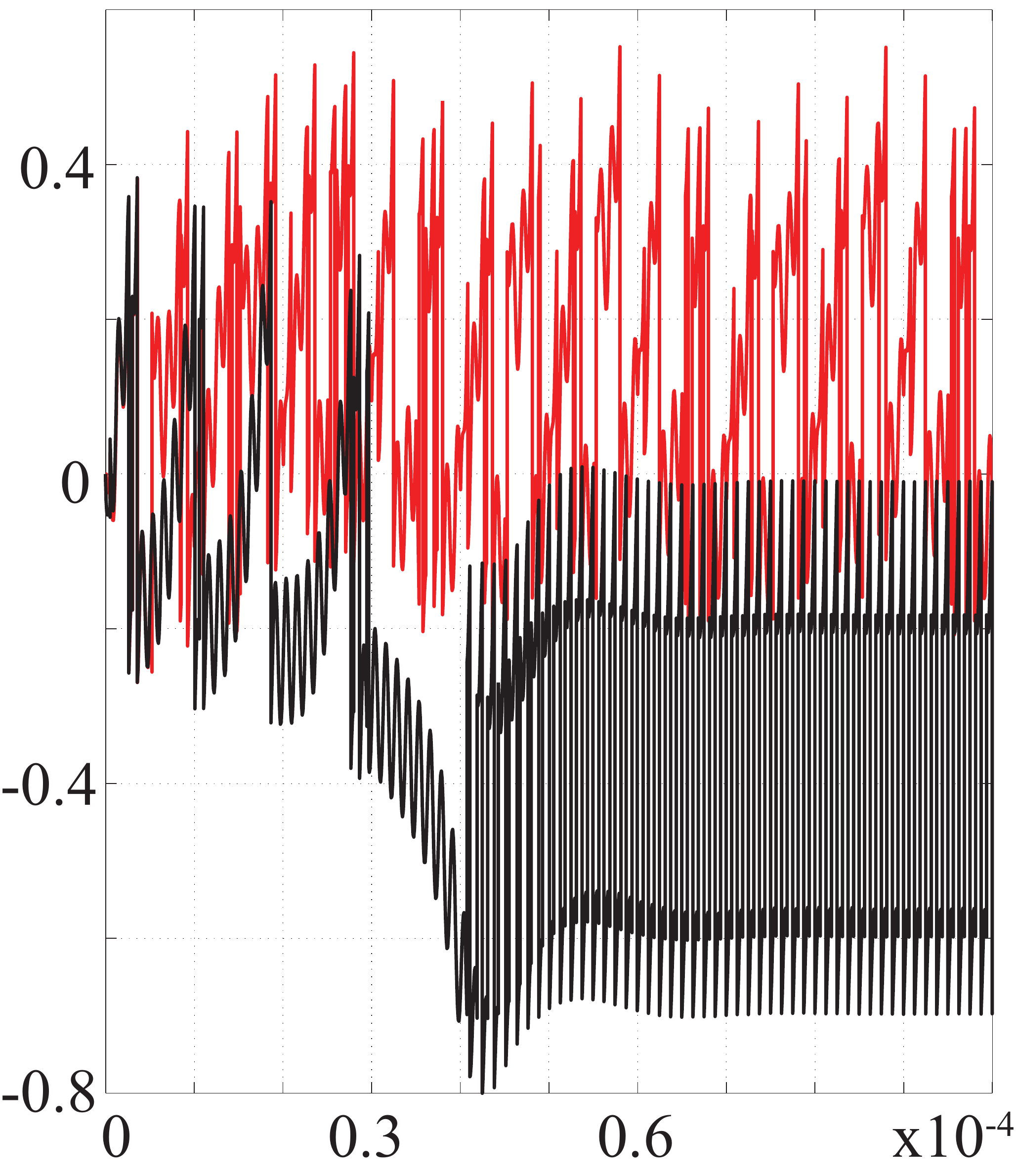}
  \caption{Loop filter output:
  VCO free-running frequency are $2.7495\cdot10^6$ rad/s;
  VCO phase shift is zero;
  initial conditions of all filters are zero;
  $\omega_{lpf} = 1.2566\cdot10^6$;
  $m_1(t) = \sign(\sin(2.7495\cdot10^6 t))$ -- red curve, $m_1(t) = 1$ -- black curve}
  \label{example_lock_unlock_with_data}
\end{figure}
\end{example}

\begin{example}
In Fig.~\ref{example_vco_phase_shift}.
is shown that initial phase of VCO may affect stability of the loop:
while physical model (black) with zero initial VCO phase
acquires lock, physical model (red)
with initial VCO phase equal to $0.8854$ rad does not acquire lock.
\begin{figure}
\centering
  \includegraphics[width=0.4\textwidth]{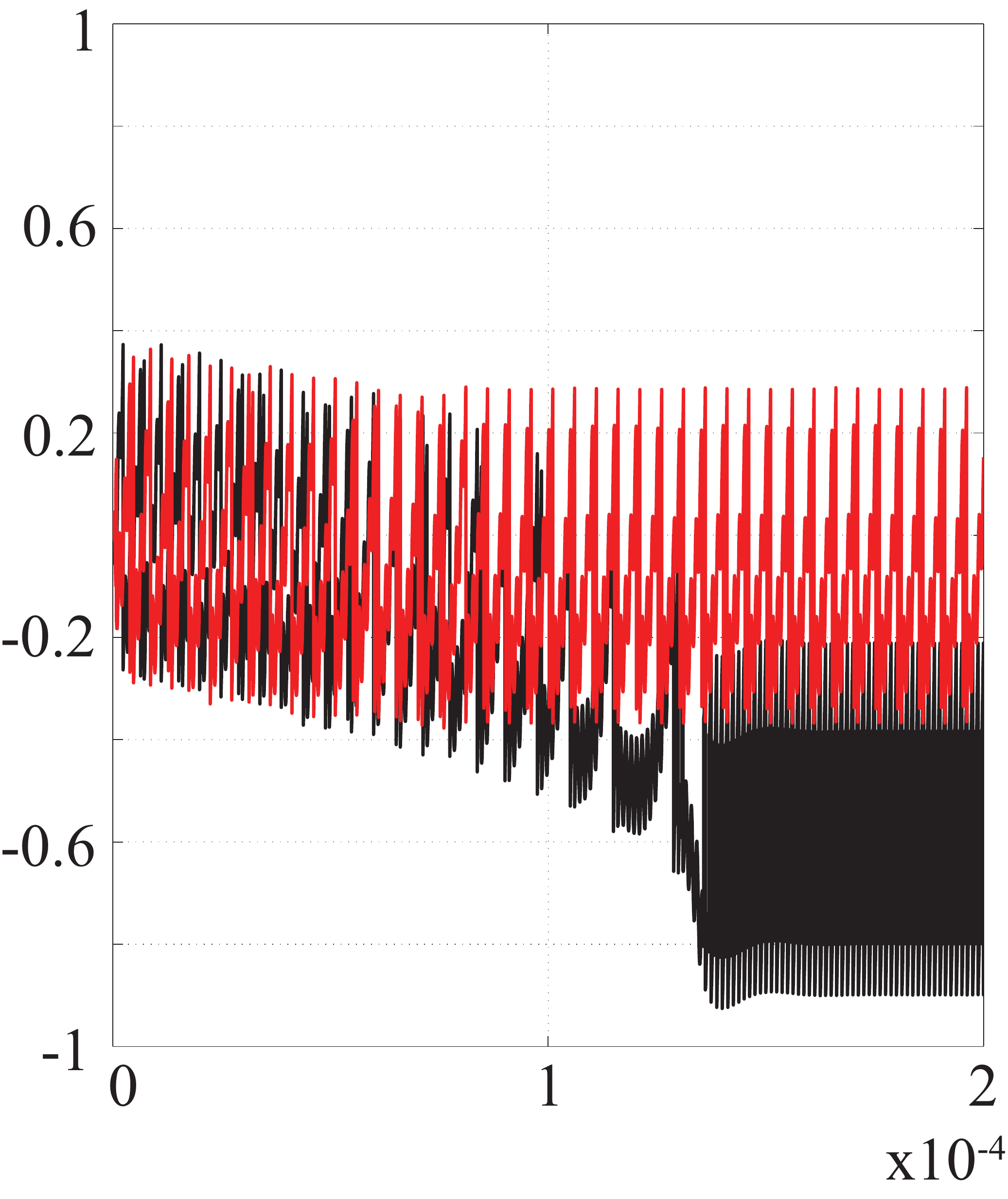}
\caption{Loop filter output for physical model:
$m_1(t) = 1$;
VCO free-running frequency is $2.8767\cdot10^6$ rad/s;
VCO phase shift is  for red curve and zero for black curve;
$\omega_{lpf} = 1.2566\cdot10^6$;
initial conditions of all filters are zero.}
\label{example_vco_phase_shift}
\end{figure}
\end{example}

\begin{example}
In Fig.~\ref{freq_def_example2}
is shown that
initial states of low-pass filters
and initial phase difference may affect stability domain:
while physical model
with low-pass filters initial states $x_1(0) = 3$, $x_2(0) = 4.2566$
and zero is out of lock (red),
classical mathematical model in the signal's phase space
with initial phase shift $\theta_{\Delta}(0)=-\frac{\pi}{4}$ rad acquires lock (black).
Therefore the consideration of classical mathematical model in the signal's phase space
(Fig.~\ref{pll-qpsk} and system \eqref{final_system})
may lead to wrong conclusion.

Since loop filter inputs corresponding to both models are equal to $1\over4$,
initial VCO control inputs (and initial VCO frequencies)
for both examples are the same.

Therefore instead of one-dimensional stability ranges
defined by $|\omega_{\Delta}|$
it is necessary to consider multi-dimensional stability domains
taking into account initial phase difference $\{\theta_{\Delta}(0)\}$.

Here VCO free-running frequency is $\omega_{\rm vco}^{free}=2.5933\cdot10^6$,
no data are being transmitted $m(t) = 1$, and initial loop filter state is zero $x(0)=\alpha_0(t) = 0$
\begin{figure}
  \includegraphics[scale=0.3]{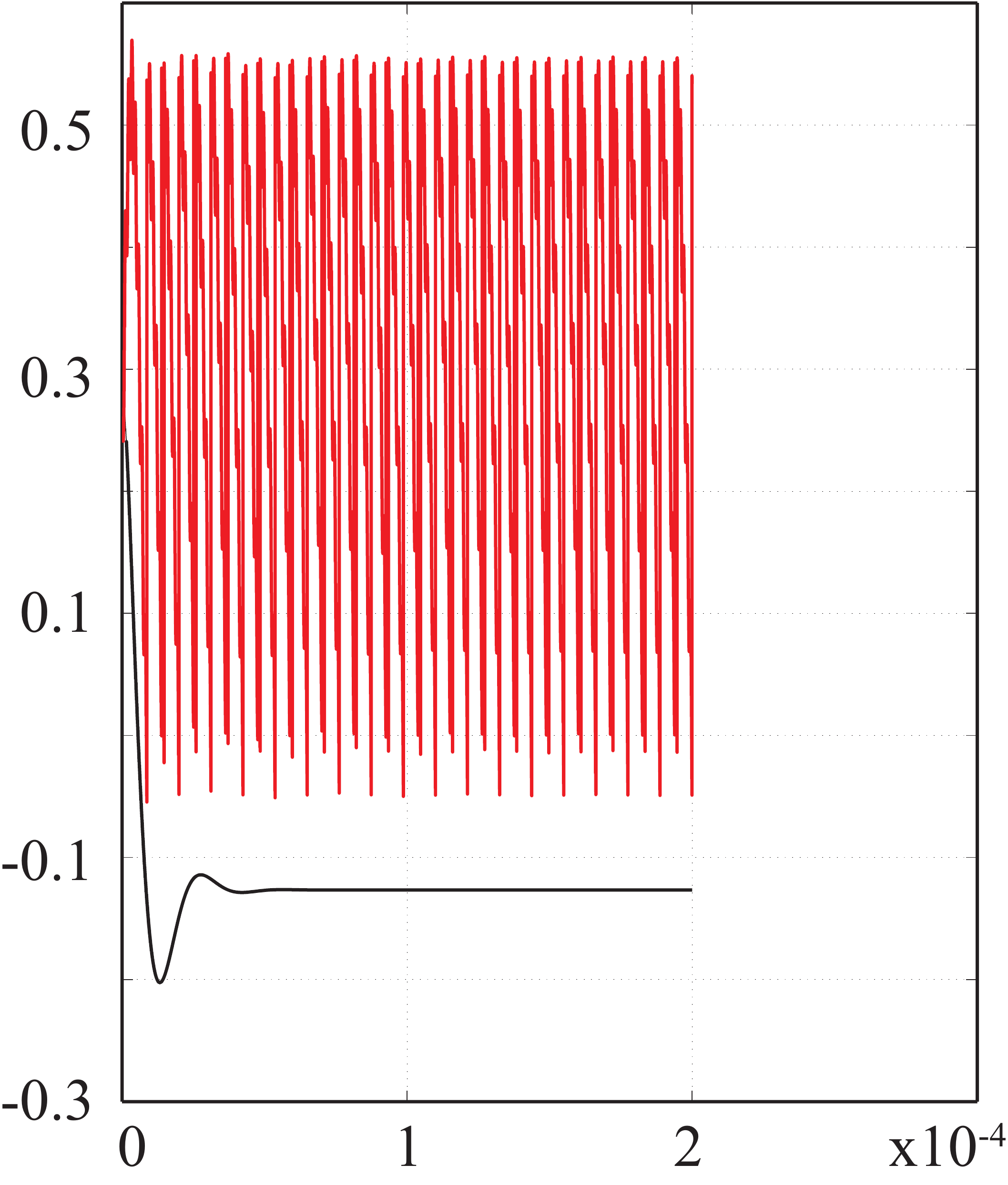}
  \caption{Loop filter output $g(t)$ for
  signal's phase space model (black curve),
  physical model (red curve).}
  \label{freq_def_example2}
\end{figure}
\end{example}

\begin{example}
In Fig.~\ref{example_lpf_phase_instablility}.
is shown that lowering corner frequency of the low-pass filter (therefore changing phase shift) may affect stability of the loop:
while signals phase model (black)
with $\omega_{lpf} = 6.2832\cdot10^5$
acquires lock,
signals phase model (red)
with $\omega_{lpf} = 1.5708\cdot10^5$ does not acquire lock.
\begin{figure}
\centering
  \includegraphics[width=0.4\textwidth]{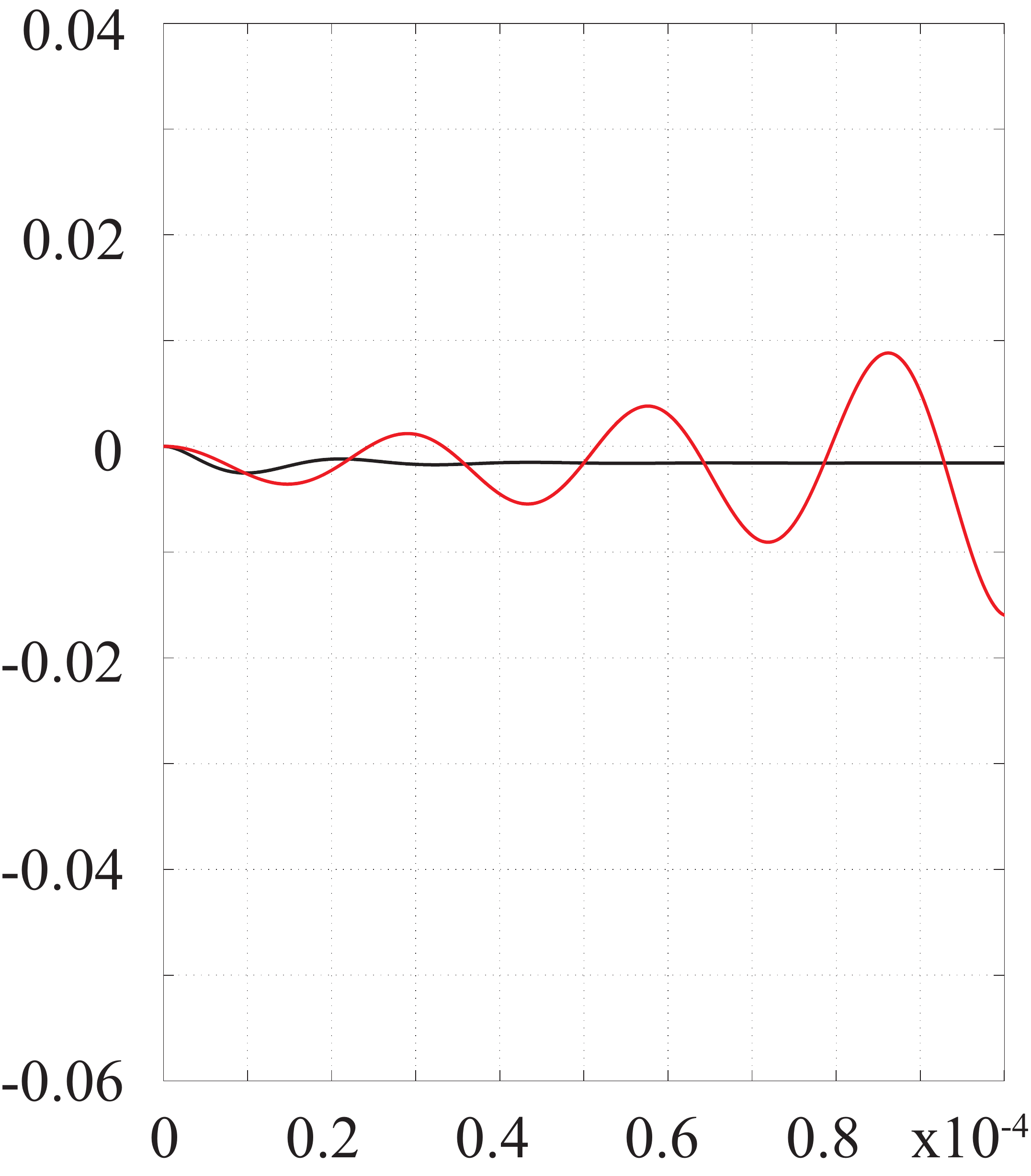}
\caption{Loop filter signal's phase model:
$m_1(t) = 1$;
VCO free-running frequency is $\omega_{\rm ref} + 1000$;
VCO phase shift is  for red curve and zero for black curve;
initial conditions of all filters are zero.}
\label{example_lpf_phase_instablility}
\end{figure}
\end{example}

\section{Conclusion}
In this survey various mathematical models of QPSK Costas loop are derived.
It is shown that the consideration of simplified mathematical models,
and the application of non rigorous methods of analysis (e.g., a simulation)
can lead to wrong conclusions concerning
the operability of \emph{physical model} of Costas loop.

\section*{\uppercase{Acknowledgements}}
 This work was supported by Russian Science Foundation (project 14-21-00041).

\bibliographystyle{ifacconf}

\end{document}